\documentclass[12pt,reqno]{amsart}
\usepackage{dsfont}
\usepackage{amssymb, amsmath, amsthm}
\usepackage{enumerate,color}
\usepackage{graphicx}
\numberwithin{equation}{section}

\newcommand{\cal}{\mathcal}
\newcommand{\bel}{\begin{equation} \label}
\newcommand{\ee}{\end{equation}}
\def\beq{\begin{equation}}
\def\eeq{\end{equation}}
\newcommand{\bea}{\begin{eqnarray}}
\newcommand{\eea}{\end{eqnarray}}
\newcommand{\beas}{\begin{eqnarray*}}
\newcommand{\eeas}{\end{eqnarray*}}

\def\Bbb R{{\rm \bf R}}
\def\C{{\rm \bf C}}
\def\proclaim#1{\vskip2mm{\bf #1}\em}
\def\endproclaim{\em \vskip2mm}
\def\tag#1{\eqno(#1)}
\def\gathered{\begin{array}{c}}
\def\endgathered{\end{array}}
\def\text{\mbox}

\title {The enclosure method for the detection of variable order in fractional diffusion equations}
\author{
Masaru IKEHATA$^\dag$
and Yavar KIAN$^*$
}

\begin{document}


\begin{abstract} 
This paper is concerned with a new type of inverse obstacle problem governed by a variable-order time-fraction diffusion equation in a bounded domain.
The unknown obstacle is a region where the space dependent variable-order of fractional time derivative 
of the governing equation deviates from a known homogeneous background one.
The observation data is given by the Neumann data of the solution of the governing equation for a specially designed Dirichlet data.
Under a suitable jump condition on the deviation, it is shown that the most recent version of {\it the time domain enclosure method} 
enables one to extract information about the geometry of the obstacle and a qualitative nature of the jump, from the observation data.

\noindent
AMS: 35R30, 35L05

\end{abstract}

\maketitle

\renewcommand{\thefootnote}{\fnsymbol{footnote}}
\footnotetext{\hspace*{-5mm} 
\begin{tabular}{@{}r@{}p{13cm}@{}} 
& Manuscript last updated: \today.
\\
$^\dag$& Laboratory of Mathematics,
Graduate School of Advanced Science and Engineering,
Hiroshima University,
Higashihiroshima 739-8527, Japan.\\
&Emeritus Professor at Gunma University, Maebashi 371-8510, Japan.\\
&E-mail: ikehataprobe@gmail.com\\
$^*$& Aix Marseille Universit\'e, Universit\'e de Toulon, CNRS, CPT, Marseille, France.
 Email: yavar.kian@univ-amu.fr\end{tabular}}

\section{Introduction}

In the present article we consider a model of anomalous diffusion described by variable order time fractional diffusion equation 
on $\Omega$  a bounded domain of $\Bbb R^3$ with $C^2$-boundary.
Namely, we fix $\alpha\in L^{\infty}(\Omega)$ satisfying $\text{ess}.\inf_{x\in\Omega}\alpha(x)>0$
and $\text{ess.}\,\sup_{x\in\Omega}\alpha (x)<1$.
Then, given $g(x,t), (x,t)\in\,\partial\Omega\times]\,0,\,\infty[$,  let
$u=u_g(x,t)$, with $(x,t)\in\,\Omega\times\,]\,0,\,\infty[$, 
denote the solution of the following 
initial boundary value problem:
$$\displaystyle
\left\{
\begin{array}{ll}
\displaystyle
(\partial_t^{\alpha(x)}-\Delta) u=0, & (x,t)\in\,\Omega\times\,]\,0,\,\infty[,\\
\\
\displaystyle
u(x,0)=0, & x\in\Omega,
\\
\\
\displaystyle
u(x,t)=g(x,t), & (x,t)\in\,\partial\Omega\times]\,0,\,\infty[.
\end{array}
\right.
\tag {1.1}
$$
Here, the symbol $\partial_t^{\alpha(x)}$
denotes the Caputo fractional derivative of order $\alpha(x)$ with respect to $t$, that is
$$\begin{array}{ll}
\displaystyle
\partial_t^{\alpha(x)}u(x,t)
=\frac{1}{\Gamma(1-\alpha(x))}
\int_0^t(t-s)^{-\alpha(x)}\partial_s u(x,s)\,ds,
& (x,t)\in\,\Omega\times\,]\,0,\,\infty[,
\end{array}
$$
where $\Gamma$ is the Gamma function.  We consider solutions $u=u_g$ of the problem   (1.1) lying in the space $C^1([0,\infty[;L^2(\Omega))\cap C([0,\infty[;H^2(\Omega)) $. 
We  denote by $\nu$ the outer unit normal vector field on $\partial\Omega$.

In this article we consider the inverse problem of determining the region of variation  and additional information about the amplitude of variation of the fractional order $\alpha(x)$ appearing in (1.1). More precisely, let $D$ be a nonempty bounded open subset of $\Omega$ with $C^2$-boundary such that $\overline D\subset\Omega$.
Assume that the order $\alpha(x)$ in (1.1) takes the form
$$\displaystyle
\alpha(x)
=\left\{
\begin{array}{ll}
\displaystyle
\alpha_0, & x\in\Omega\setminus D,\\
\\
\displaystyle
\alpha_0+h(x), & x\in D,
\end{array}
\right.
\tag {1.2}
$$
where $\alpha_0\in\,]0,\,1[$ and
the function $h$ belongs to $L^{\infty}(D)$ and satisfies
$$
\displaystyle
-\alpha_0<\text{ess.inf}_{x\in D}\,h(x)\le \text{ess.sup}_{x\in D}\,h(x)<1-\alpha_0.
$$

We impose the jump condition (A.I)/(A.II) of $\alpha$ from $\alpha_0$ across $\partial D$ as follows.
$$\begin{array}{lllll}
\displaystyle
\text{(A.I)} &   \exists C>0 & \exists \gamma\ge 0 &  \text{$h(x)\ge C\,\text{dist}\,(x,\partial D)^{\gamma}$ a.e. $x\in D$;}
\\
\\
\displaystyle
\text{(A.II)} &  \exists C>0 & \exists \gamma\ge 0 &  \text{$-h(x)\ge C\,\text{dist}\,(x,\partial D)^{\gamma}$ a.e. $x\in D$.}
\end{array}
$$
To briefly describe the difference between the two conditions sometimes we write $\alpha>>\alpha_0$ if (A.I) is satisfied; $\alpha<<\alpha_0$ if (A.II) is satisfied. Our inverse problem can be stated as follows.

$\quad$

{\bf\noindent Problem.}  Assume that $\alpha_0$ is known and both $D$ and $h$ are unknown.
Given $g$ (to be specified later) we extract information about the
location and shape of $D$ and qualitative property of $h$ from the Neumann data $\partial_\nu u_g$
on $\partial\Omega$ over the time interval $]0,\,\infty[$.

$\quad$

Recall that the initial boundary value problem (1.1) is frequently used as a model for anomalous diffusion in complex media
with  applications in different fields such as   geophysics, environmental 
 and biological problems. Such diffusion process are often described by problem (1.1) with a constant order $\alpha$ (see  \cite{AG,CSLG}).
 However, in some complex media the presence of heterogeneous regions displays space inhomogeneous variations and the constant order fractional dynamic models are not robust for long times (see \cite{FS,FH}). For such problems, the variable order time-fractional model is considered as  more relevant for describing the space-dependent anomalous diffusion process (see e.g. \cite{SCC}). Indeed, several variable order diffusion models have been successfully applied in different problems of sciences and engineering, including Chemistry \cite{CZZ}, 
Rheology \cite{SdV}, Biology \cite{GN}, Hydrogeology \cite{AON} and Physics \cite{SS, ZLL}. In this context, the goal of our inverse problem is to determine information about the variable order $\alpha$ which play a fundamental role in the anomalous mechanism leading to the model (1.1).

The inverse problems of determining  fractional orders,   which is one of the
most important inverse problems for fractional diffusion equations,  have been extensively studied these last decades. We refer to \cite{LLY19} for a survey about this topic (see also \cite{JR} for an overview of inverse problems for fractional diffusion equations). Without being exhaustive we can mention the works of \cite{AA20,CNY,HNWY,Ja,JaK,K1,LW,LIY,LY,LZ,Y21} devoted to the determination of single or multiple constant fractional orders, sometimes together with other parameters (coefficients or internal sources),  from several class of observational data. We mention also the recent works  \cite{JK1,JK2} where the determination of constant fractional order  have been studied in the context of an unknown medium (unknown source, coefficients, domain...).
All the above mentioned results have been devoted to the determination of constant  fractional order. The only result that we are aware of dealing with the determination of variable fractional order depending on the space variable can be found in \cite{KSY}. Here the authors proved the determination of general order $\alpha\in L^\infty(\Omega)$  from the knowledge of Neumann data $\partial_\nu u_g$ on $\partial\Omega\times\,]0,\,T[$ with an arbitrary fixed $T>0$
for {\it infinitely many} input $g$ having the form $g(x,t)=t^{\kappa}g(x)$ with a constant $\kappa\in ]2,\,\infty[$.  
The aim of the present article is to prove the detection of the region of variation and the amplitude of variation of $\alpha$   by using the {\it enclosure method} initiated in \cite{E00} where the infinite boundary measurements under consideration in \cite{KSY} are replaced by a single boundary measurement for some class of suitable input $g$. We apply here the most recent version of \textit{the time domain enclosure method} developed in \cite{IEO4,IE06}.

\subsection{Statement of the main result}

Now let us describe the main result of this paper. For this purpose, we start by introducing the class of input under consideration in (1.1).
Let $\eta>0$ and $0<R_1<R_2$.  Let $B_{\eta}$ be the open ball with radius $\eta$ such that $\overline{B_{\eta}}\cap\overline\Omega=\emptyset$.
Let $B_{R_1}$ and $B_{R_2}$ be two concentric balls with radius $R_1$ and $R_2$, respectively such that $\Omega\subset B_{R_1}$ (see Figure 1).
Let $m=0,1,\cdots$.
For a complex number $\tau$ with $\text{Re}\,\tau>1$, we choose the solution $w_{\star,m}^0(x)=w_{\star,m}^0(x,\tau)\in H^2(\Bbb R^3)$ with $\star=\text{ext}, \text{int}$ of the equation
$$\begin{array}{ll}
\displaystyle
(\Delta-\tau^{\alpha_0}\,)w_{\star,m}^0+\tau^{\alpha_0-1}\,\Psi_{\star,m}(x)=0, & \displaystyle x\in\Bbb R^3,
\end{array}
$$
where
$$\left\{
\begin{array}{l}
\displaystyle
\Psi_{\text{ext},m}(x)=(\eta^2-\vert x-p\vert^2)^m\chi_{B_{\eta}}(x),\\
\\
\displaystyle
\Psi_{\text{int},m}(x)=(R_2^2-\vert x-p\vert^2)^m(\vert x-p\vert^2-R_1^2)^m\,\chi_{B_{R_2}\setminus B_{R_1}}(x),
\end{array}
\right.
$$
the point $p$ denotes the center of $B_{\eta}$ when $\star=\text{ext}$; the common center of $B_{R_1}$ and $B_{R_2}$
when $\star =\text{int}$.  Note that both $\Psi_{\text{ext},m}(x)$ and $\Psi_{\text{int},m}(x)$ are non-negative for all $x\in\Bbb R^3$.
Note that the restriction of $w^0_{\star,m}$ onto $\Omega$ satisfies
$$\begin{array}{ll}
\displaystyle
(\Delta-\tau^{\alpha_0})w^0_{\star,m}=0, & x\in\Omega.
\end{array}
\tag {1.3}
$$
In the present article we consider the following class of input $g$
$$\begin{array}{lll}
\displaystyle
g_{\star,m}(x,t)=\frac{e^t}{2\pi}\,
\int_{-\infty}^{\infty}
e^{its}(1+is)^{-5}w_{\star,m}^0(x,1+is)ds, & x\in\partial\Omega, & t\in[0,\,\infty[,
\end{array}
\tag {1.4}
$$
where $\star=\text{ext},\,\text{int}$.

In order to state our main result we need to consider first the forward problem. Namely,  we consider solutions of  (1.1) lying in the space $C^1([0,\infty[;L^2(\Omega))\cap C([0,\infty[;H^2(\Omega)) $. This means that we consider solutions of  (1.1) in a strong sense as stated in \cite[Definition 2.2]{KY}. 
In addition to this property, we will show in the next Proposition 1.2 that $e^{-t}g_{\star,m}\in  L^{\infty}(\Bbb R_+;H^{\frac{3}{2}}(\partial\Omega))$ and,
for all $\tau>1$, $e^{-\tau t}u\in L^1(\Bbb R_+; L^2(\Omega))$. Moreover, we will show that, for all $\tau\in \C$ satisfying $\text{Re}\tau>1$, the Laplace transform 
$$\displaystyle
\hat{u}(\,\cdot\,,\tau)=\int_0^{\infty}e^{-\tau t}u(t,\,\cdot\,)\,dt
$$
of $u$ is lying in $H^2(\Omega)$ and it solves the boundary value problem
$$\left\{
\begin{array}{ll}
\displaystyle
(\Delta-\tau^{\alpha(x)})\hat{u}(x,\tau)=0, & x\in\Omega,\\
\\
\displaystyle
\hat{u}(x,\tau)=\widehat{g_{\star,m}}(x,\tau), & x\in\partial\Omega.
\end{array}
\right.
\tag {1.5}
$$
All these properties are stated in  the following proposition.
\proclaim{\noindent Proposition 1.1.}
Let $g_{\star,m}$ be given by (1.4).  Then, we have $e^{-t}g_{\star,m}\in L^{\infty}(\Bbb R_+; H^{\frac{3}{2}}(\partial\Omega))$
and, for all complex  $\tau$ satisfying $\textrm{Re}\tau>1$, we have
$$\begin{array}{ll}
\displaystyle
\widehat{g_{\star,m}}(x,\tau)=\tau^{-5}w^0_{\star,m}(x,\tau), & x\in\partial\Omega.
\end{array}
\tag {1.6}
$$
Moreover, for $g=g_{\star,m}$, the problem (1.1) admits a unique strong solution\\ 
$u_{\star,m}\in C^1([0,\infty[;L^2(\Omega))\cap C([0,\infty[;H^2(\Omega)) $  satisfying  $e^{- t}u\in W^{1,\infty}(\Bbb R_+;H^2(\Omega))$. Finally, for all complex  $\tau$ satisfying $\textrm{Re}\tau>1$, the Laplace transform in time $\widehat{u_{\star,m}}$ of $u_{\star,m}$ solves (1.5), we have $e^{-\tau t}\partial_{\nu}u_{\star,m}\in L^1(\Bbb R_+;H^{\frac{3}{2}}(\partial\Omega))$ and
$$\begin{array}{ll}
\displaystyle
\widehat{\partial_{\nu}u_{\star,m}}(x,\tau)=\partial_{\nu}\widehat{u_{\star,m}}(x,\tau), & x\in\partial\Omega.
\end{array}
\tag {1.7}
$$

\endproclaim

We mention that the only other works that we are aware of dealing with the existence of solutions of (1.1) with variable order can be found in \cite{K2,KSY} and only \cite{KSY} considered this problem with non-homogenous boundary condition. In Proposition 1.1,  we extend the analysis of \cite{KSY} to more general class of Dirichlet boundary conditions of the form (1.4) and, in contrast to \cite{KSY} who considered solutions defined in terms of Laplace transform in time, we prove the unique existence of strong solutions of (1.1).

Applying Proposition 1.1,  we introduce the following indicator function which is one of the key ingredient of the {\it enclosure method}.

$\quad$

{\bf\noindent Definition 1.2.}
Define
$$\begin{array}{ll}
\displaystyle
I_{\star,m}(\tau)=\int_{\partial\Omega}
\left(\widehat{\partial_{\nu}u_{\star,m}}(x,\tau)-\tau^{-5}\partial_{\nu}w^0_{\star,m}(x,\tau)\,\right)
\tau^5w^0_{\star,m}(x,\tau)\,dS(x), & \tau>1.
\end{array}
\tag {1.8}
$$

$\quad$

This indicator function can be computed from the data $\partial_\nu u_g$ on $\partial\Omega$ over
time interval $]0,\,\infty[$ which is the Neumann data of the  solution of (1.1) with $g=g_{\star,m}$.
Applying Proposition 1.1, we can transform the indicator function (1.8) in the following way.

Let the function $w_{\star,m}$ belongs to $H^2(\Omega)$ and satisfies
$$\left\{\begin{array}{ll}
\displaystyle
(\Delta-\tau^{\alpha(x)})w_{\star,m}=0,& x\in\Omega,\\
\\
\displaystyle
w_{\star,m}(x,\tau)=w_{\star,m}^0(x,\tau), & x\in\partial\Omega.
\end{array}
\right.
\tag {1.9}
$$
In view of (1.6), we have
$$\begin{array}{lll}
\displaystyle
w^0_{\star,m}(x,\tau)=\tau^5\widehat{g_{\star,m}}(x,\tau), & \tau>1, & x\in\partial\Omega
\end{array}
$$
and (1.5) together with (1.9) yields
$$\begin{array}{lll}
\displaystyle
w_{\star,m}(x,\tau)=\tau^5\widehat{u_{\star,m}}(x,\tau), & \tau>1, & x\in\Omega.
\end{array}
$$
Applying this together with (1.7) to (1.8), we obtain the more familiar expression of the indicator function
$$\begin{array}{ll}
\displaystyle
I_{\star,m}(\tau)
=\int_{\partial\Omega}\,(\partial_{\nu}w_{\star,m}(x,\tau)-\partial_\nu w^0_{\star,m}(x,\tau))w^0_{\star,m}(x,\tau)\,dS(x),
& \tau>1.
\end{array}
\tag {1.10}
$$

Fixing $K_{\star}=\text{supp}\,\Psi_{\star,m}$, we get
$$\displaystyle
K_{\star}
=
\left\{
\begin{array}{ll}
\overline{B_{\eta}}, & \text{if $\star=\text{ext}$,}
\\
\\
\displaystyle
\overline{B_{R_2}}\setminus B_{R_1}, & \text{if $\star=\text{int}$.}
\end{array}
\right.
$$
Besides we have
$$\displaystyle
\text{dist}\,(K_{\star},D)=
\left\{
\begin{array}{ll}
\displaystyle
\text{dist}\,(p,D)-\eta, & \text{if $\star=\text{ext}$,}
\\
\\
\displaystyle
R_1-R_D(p), & \text{if $\star=\text{int}$,}
\end{array}
\right.
$$
where $\text{dist}(p,D)=\inf_{x\in D}\vert x-p\vert$ and $R_D(p)=\sup_{x\in D}\,\vert x-p\vert$.
Thus knowing the value of $\text{dist}\,(K_{\star},D)$ is equivalent to that of $\text{dist}\,(p,D)$/$R_D(p)$ if $\star=\text{ext}$/$\text{int}$.
Note that the sphere $\vert x-p\vert=\text{dist}\,(p,D)$ is the largest one whose exterior contains $D$:
the sphere $\vert x-p\vert=R_D(p)$ is the smallest one whose interior contains $D$ (see Figure 1 for more detail).

\begin{figure}[!ht]
  \centering
  \includegraphics[width=0.4\textwidth]{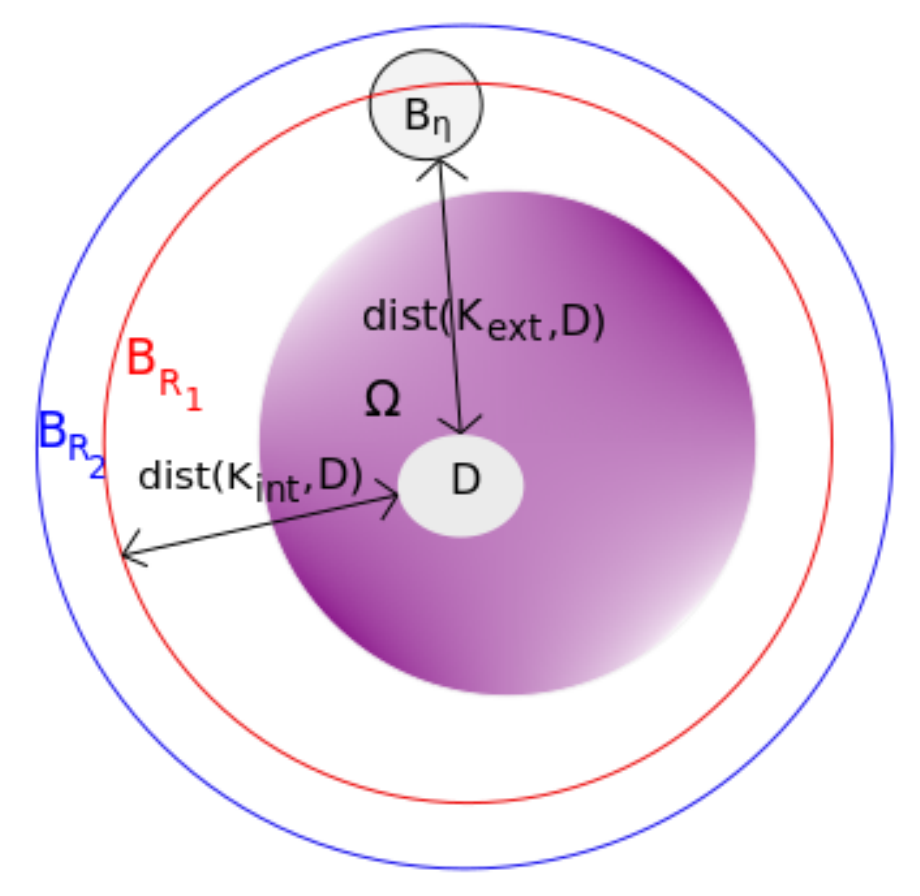}
  \caption{The sets $\Omega$, $D$ and $\text{dist}\,(K_{\star},D)$.  }
  \label{fig1}
\end{figure}

Now we state the main result of this paper.

\proclaim{\noindent Theorem 1.1}  Let $T$ be an arbitrary positive number.
We have
$$\begin{array}{l}
\displaystyle
\lim_{\tau\rightarrow\infty}e^{\tau^{\frac{\alpha_0}{2}}T}I_{\star,m}(\tau)
=
\left\{
\begin{array}{ll}
\displaystyle
0 & \text{if $T<2\,\text{dist}\,(K_{\star},D)$,}
\\
\\
\displaystyle
\infty & \text{if $\alpha>>\alpha_0$ and $T>2\text{dist}\,(K_{\star},D)$,}
\\
\\
\displaystyle
-\infty & \text{if $\alpha<<\alpha_0$ and $T>2\text{dist}\,(K_{\star},D)$.}
\end{array}
\right.
\end{array}
\tag {1.11}
$$
If (A.I) or (A.II) are satisfied, then there exists a positive number $\tau_0$ such that, for all $\tau\ge\tau_0$
$\vert I_{\star,m}(\tau)\vert>0$ and we have the one line formula
$$\displaystyle
\lim_{\tau\longrightarrow\infty}
\tau^{-\frac{\alpha_0}{2}}
\log \vert I_{\star,m}(\tau)\vert=-2\text{dist}\,(K_{\star},D).
\tag {1.12}
$$

\endproclaim

Let us observe that  Theorem 1.1 give several important information about the domain of variation $D$ and the amplitude of variation $h$ of the variable order $\alpha$.
Namely, formula (1.11) gives a target distinction and range estimate at the same time, that means one can distinguish whether $\alpha>>\alpha_0$ or $\alpha<<\alpha_0$
together with $T>2\,\text{dist}\,(K_{\star},D)$ or $T<2\,\text{dist}\,(K_{\star},D)$ (see figure 1 for more detail)
by using the asymptotic behavior of the indicator function as $\tau\rightarrow\infty$.
Formula (1.12) gives us a direct way of extracting information about the geometry of $D$ from the indicator function.

To the best of our knowledge, in Theorem 1.1 we obtain the first result of extraction of information about 
the variable order $\alpha$ from a single boundary measurement of the solution of (1.1). 
Indeed, the only other work treating this type of problem can be found in \cite{KSY} where the authors considered the problem 
of recovering the full knowledge of $\alpha$ itself from infinite boundary measurements. 
We give in this article an application of the \textit{enclosure method}, considered so far mainly for inverse source or inverse obstacle problem 
\cite{Isource, Icavity}, to a new class of inverse obstacle problem, that is the problem of extracting information about the region of the jump
of space dependent variable order of fractional time derivative 
in the governing equation from the background one.

This article is organized as follows. Section 2 is devoted to the proof of the result about the froward problem stated in Proposition 1.1 where we show the unique existence of strong solutions of (1.1) having some specific properties. In Section 3, we complete the proof of our main result stated in Theorem 1.1 by assuming Lemma 3.2 whose proof is postponed to Section 4. Finally, in Section 5 we give some additional remarks about our results with possible extension of our analysis.

\section{Proof of Proposition 1.1}

In all this proof $C>0$ will be a constant independent of $\tau$ that may change from line to line.

\subsection{Proof of (1.6)}

Let us first observe that for all complex $\tau\in\C_+:=\{z\in \C:\ \text{Re}\,z>0\}$ the Fourier transform ${\cal F}w^0_{\star,m}$ in $x$ of $w^0_{\star,m}$ is given by
$$\displaystyle
{\cal F}w^0_{\star,m}(\xi,\tau)
=\frac{\tau^{\alpha_0-1}\,{\cal F}\Psi_{\star,m}(\xi)}{\vert\xi\vert^2+\tau^{\alpha_0}}.
$$
We fix $r=\vert 1+\tau\vert>1$ and $\theta\in\,]-\frac{\pi}{2},\,\frac{\pi}{2}[$ such that $1+\tau=re^{i\theta}$.
Using the fact that $\Psi_{\star,m}\in\,L^2(\Bbb R^3)$, we get
$$\begin{array}{ll}
\displaystyle
\Vert w^0_{\star,m}(\,\cdot\,,1+\tau)\Vert^2_{H^2(\Bbb R^3)}
&
\displaystyle
=
\vert 1+\tau\vert^{2(\alpha_0-1)}
\,\int_{\Bbb R^3}\left(\frac{\vert\xi\vert^2+1}{\vert\vert\xi\vert^2+(1+\tau)^{\alpha_0}\vert}\,\right)^2\,\vert{\cal F}\Psi_{\star,m}(\xi)\vert^2\,d\xi
\\
\\
\displaystyle
&
\displaystyle
\le
\vert 1+\tau\vert^{2(\alpha_0-1)}
\,\int_{\Bbb R^3}\left(\frac{\vert\xi\vert^2+1}{\vert\xi\vert^2+r^{\alpha_0}\cos(\alpha_0\theta)}\,\right)^2\,\vert{\cal F}\Psi_{\star,m}(\xi)\vert^2\,d\xi
\\
\\
\displaystyle
&
\displaystyle
\le
\vert 1+\tau\vert^{2(\alpha_0-1)}
\,\int_{\Bbb R^3}\left(\frac{\vert\xi\vert^2+1}{\vert\xi\vert^2+\cos(\frac{\alpha_0\pi}{2})}\,\right)^2\,\vert{\cal F}\Psi_{\star,m}(\xi)\vert^2\,d\xi
\\
\\
\displaystyle
&
\displaystyle
\le
C\vert 1+\tau\vert^{2(\alpha_0-1)}
\,\int_{\Bbb R^3}\,\vert{\cal F}\Psi_{\star,m}(\xi)\vert^2\,d\xi
\\
\\
\displaystyle
&
\displaystyle
\le
C\vert 1+\tau\vert^{2(\alpha_0-1)}
\Vert\Psi_{\star,m}\Vert^2_{L^2(\Bbb R^3)}.
\end{array}
$$
Here we have used the fact that $\cos(\frac{\alpha_0\pi}{2})>0$ since $\alpha_0\in\,]0,\,1[$.
Thus,  we find
$$\begin{array}{ll}
\displaystyle
\Vert(1+\tau)^{-5}w^0_{\star,m}(\,\cdot\,,1+\tau)\vert_{\partial\Omega}\Vert_{H^{\frac{3}{2}}(\partial\Omega)}
&
\displaystyle
\le C\vert 1+\tau\vert^{-5}\Vert w^0_{\star,m}(\,\cdot\,,1+\tau)\Vert_{H^2(\Bbb R^3)}
\\
\\
\displaystyle
&
\displaystyle
\le
C\vert 1+\tau\vert^{\alpha_0-6},\quad \tau\in\C_+.
\end{array}
\tag {2.1}
$$
This together with (1.4) yields $e^{-t}g_{\star,m}\in L^{\infty}(\Bbb R_+;H^{\frac{3}{2}}(\partial\Omega))$.
Moreover, using similar arguments as above, one can check that the map $\tau\mapsto w^0_{\star,m}(\,\cdot\,,1+\tau)$ 
is holomorphic with respect to $\tau\in\C_+$ as a map taking values in $H^2(\Bbb R^3)$.
And it follows that the map $\tau\mapsto (1+\tau)^{-5}\,w^0_{\star,m}(\,\cdot\,,1+\tau)\vert_{\partial\Omega}$
is holomorphic with respect to $\tau\in\C_+$ as a map taking values in $H^{\frac{3}{2}}(\partial\Omega)$.
Therefore, applying \cite{R} Theorem 19.2 and the note, we deduce that
$$\begin{array}{lll}
\widehat{g_{\star,m}}(x,1+\tau)=\widehat{e^{-t}g_{\star,m}}(x,\tau)=(1+\tau)^{-5}w^0_{\star,m}(x,1+\tau), & x\in\partial\Omega,\ \tau\in\C_+.
\end{array}
$$
This identity clearly implies (1.6).

\subsection{Proof of the unique existence of  solutions of (1.1) with $g=g_{\star,m}$ and (1.7)}

For $\tau\in\C_+$ let $v(\,\cdot\,,\tau)$ be the solution of
$$\left\{
\begin{array}{ll}
\displaystyle
\Delta v(x,\tau)-\tau^{\alpha(x)}\,v(x,\tau)=0, & x\in\Omega,\\
\\
\displaystyle
v(x,\tau)=\widehat{g_{\star,m}}(x,\tau), & x\in\partial\Omega.
\end{array}
\right.
\tag {2.2}
$$
One can split $v(x,1+\tau)$ into two terms 
$$\begin{array}{ll}
\displaystyle
v(x,1+\tau)=(1+\tau)^{-5}\,w^0_{\star,m}(x,1+\tau)+z(x,1+\tau), & x\in\Omega,\ \tau\in\C_+
\end{array}
$$
where $z(\,\cdot\,,1+\tau)$ solves
$$\left\{
\begin{array}{ll}
\displaystyle
\Delta z(x,1+\tau)-(1+\tau)^{\alpha(x)}\,z(x,1+\tau)=-(1+\tau)^{-5}G(x,1+\tau), & x\in\Omega,\\
\\
\displaystyle
z(x,1+\tau)=0, & x\in\partial\Omega
\end{array}.
\right.
$$
and 
$$\begin{array}{ll}
\displaystyle
G(x,1+\tau)=\Delta w^0_{\star,m}(x,1+\tau)-(1+\tau)^{\alpha(x)}w^0_{\star,m}(x,1+\tau), & x\in\Omega.
\end{array}
$$
Therefore, fixing $\Delta_0$ the Laplacian with Dirichlet boundary condition acting $L^2(\Omega)$, for all  $\tau\in\C_+$,
we obtain
$$\displaystyle
z(\,\cdot\,,1+\tau)
=(1+\tau)^{-5}
(-\Delta_0+(1+\tau)^{\alpha(x)})^{-1}[\Delta w^0_{\star,m}(\,\cdot\,,1+\tau)-(1+\tau)^{\alpha(x)}w^0_{\star,m}(1+\tau,\,\cdot\,)].
$$
Following Proposition 2.1 in \cite{KSY},
we have
$$\displaystyle
\Vert(-\Delta_0+(1+\tau)^{\alpha(x)})^{-1}\Vert_{{\cal B}(L^2(\Omega))}\le C\vert 1+\tau\vert^{\alpha^0_U-2\alpha^0_L},
$$
where $0<\alpha_L^0=\text{ess.inf}_{x\in\Omega}\alpha(x)\le\text{ess.sup}_{x\in\Omega}\alpha(x)=\alpha^0_U<1$.
It follows that, for all $v\in L^2(\Omega)$, we have
$$\begin{array}{ll}
\displaystyle
\Vert(-\Delta_0+(1+\tau)^{\alpha(x)})^{-1}v\Vert_{H^2(\Omega)}
&
\displaystyle
\le
C\Vert-\Delta_0(-\Delta_0+(1+\tau)^{\alpha(x)})^{-1}v\Vert_{L^2(\Omega)}
\\
\\
\displaystyle
&
\displaystyle
\le
C\Vert v-(1+\tau)^{\alpha(x)}(-\Delta_0+(1+\tau)^{\alpha(x)})^{-1}v\Vert_{L^2(\Omega)}
\\
\\
\displaystyle
&
\displaystyle
\le
C\left(1+\vert 1+\tau\vert^{\alpha^0_U}\Vert(-\Delta_0+(1+\tau)^{\alpha(x)})^{-1}\Vert_{{\cal B}(L^2(\Omega))}\right)\,\Vert v\Vert_{L^2(\Omega)}
\\
\\
\displaystyle
&
\displaystyle
\le C\vert 1+\tau\vert^{2\alpha^0_U-2\alpha^0_L}\,\Vert v\Vert_{L^2(\Omega)}.
\end{array}
$$
From this estimate and (2.1), we obtain
$$\begin{array}{ll}
\displaystyle
\Vert z(\,\cdot\,,1+\tau)\Vert_{H^2(\Omega)}
&
\displaystyle
\le C\vert 1+\tau\vert^{2\alpha^0_U-2\alpha^0_L-5}
\left(\Vert w^0_{\star,m}(\,\cdot\,,1+\tau)\Vert_{H^2(\Bbb R^3)}
\right.
\\
\\
\displaystyle
&
\displaystyle
\,\,\,
\left.
+\vert 1+\tau\vert^{\alpha^0_U}\Vert w^0_{\star,m}(\,\cdot\,,1+\tau)\Vert_{L^2(\Bbb R^3)}
\right)
\\
\\
\displaystyle
&
\displaystyle
\le
C\vert 1+\tau\vert^{3\alpha^0_U+\alpha_0-6}.
\end{array}
$$
This estimate implies that the solution $v(\,\cdot\,,1+\tau)$ of (2.2) satisfies the following estimate
$$\displaystyle
\Vert v(\,\cdot\,,1+\tau)\Vert_{H^2(\Omega)}\le C\vert 1+\tau\vert^{3\alpha^0_U+\alpha_0-6},
\tag {2.3}
$$
where we recall that $3\alpha^0_U+\alpha_0-6<-2$.
Besides, in view of Proposition 2.1 in \cite{KSY}, we know that the map $\tau\mapsto (-\Delta_0+(1+\tau)^{\alpha(x)})^{-1}$
is holomorphic with respect to $\tau\in\C_+$ as a map taking values in ${\cal B}(L^2(\Omega))$ and
from the above properties we deduce that $\tau\mapsto v(\,\cdot\,,1+\tau)$ is holomorphic with respect to $\tau\in\C_+$ 
as a map taking values in $H^2(\Omega)$.

Thus, applying again Theorem 19.2 and note in \cite{R} as well as estimate (2.3), we deduce that for
$$\begin{array}{lll}
\displaystyle
u_{\star,m}(x,t)=\frac{e^t}{2\pi}\,\int_{-\infty}^{\infty}\,e^{its}v(x, 1+is)\,ds, & x\in\Omega, & t\in\,]0,\,\infty[,
\end{array}
$$
we have $e^{-t}u_{\star,m}\in\,W^{1,\infty}(\Bbb R_+;H^2(\Omega))\cap C^1([0,\infty[;L^2(\Omega))\cap C([0,\infty[;H^2(\Omega))$ and
$$\begin{array}{lll}
\displaystyle
\widehat{u_{\star,m}}(x,1+\tau)=\widehat{e^{-t}u_{\star,m}}(x,\tau)=v(x,1+\tau), & \tau\in\C_+,\,x\in\Omega.
\end{array}
$$
This proves that for all complex $\tau$ satisfying $\text{Re}\,\tau>1$,  $\widehat{u_{\star,m}}(\cdot,\tau)$ solves the boundary value problem
(1.5). Therefore, in order to complete the proof of the proposition we need to prove that $u_{\star,m}$ is the unique solution of (1.1) for $g=g_{\star,m}$ satisfying $e^{- t}u_{\star,m}\in W^{1,\infty}(\Bbb R_+; H^2(\Omega))$.  Let us show first that $u_{\star,m}$ solves (1.1)  for $g=g_{\star,m}$. Note  that, in view of estimate (1.7), we have
$$u_{\star,m}(x,0)=\frac{1}{2\pi}\,\int_{-\infty}^{\infty}\,v(x, 1+is)\,ds,\quad x\in\Omega.$$
On the other hand, for all $R>0$, fixing  the contour  $C_R=\{1+Re^{i\theta}:\ \theta\in[-\pi/2,\pi/ 2]\}$ and applying the residue theorem we deduce that
$$\frac{1}{2\pi}\int_{-R}^{R}\,v(x, 1+is)\,ds=\frac{1}{2i\pi}\int_{ C_R}\,v(x, 1+\tau)\,d\tau\quad x\in\Omega.$$
Now sending $R\to\infty$ and applying estimate (2.3), we get
$$\begin{array}{lll}\|u_{\star,m}(\cdot,0)\|_{L^2(\Omega)}&
\displaystyle
=\left\|\frac{1}{2\pi}\,\int_{-\infty}^{\infty}\,v(\cdot, 1+is)\,ds\right\|_{L^2(\Omega)}\\
\\
&\displaystyle\leq \limsup_{R\to\infty}\left\|\frac{1}{2i\pi}\int_{ C_R}\,v(\cdot, 1+\tau)\,d\tau\right\|_{L^2(\Omega)}\\
\\
\ &\displaystyle\leq\limsup_{R\to\infty}\frac{1}{2\pi}\int_{-\pi/2}^{\pi/2}R\|v(\cdot, 1+Re^{i\theta})\|_{L^2(\Omega)}d\theta=0.\end{array}$$
It follows that $u_{\star,m}(\cdot,0)=0$. Using the fact that $e^{-t}u_{\star,m}\in\,W^{1,\infty}(\Bbb R_+;H^2(\Omega))$ and applying the properties of fractional derivative (see e.g. \cite[pp. 80]{P}), for all complex $\tau$ satisfying $\text{Re}\,\tau>1$ we deduce that 
$$\tau^\alpha \widehat{u_{\star,m}}(\cdot,\tau)-\tau^{\alpha-1}u_{\star,m}(\cdot,0)=\widehat{\partial_t^\alpha u_{\star,m}}(\cdot,\tau).$$
Combining this with the fact that $u_{\star,m}(\cdot,0)=0$, we deduce that for $y=\partial_t^\alpha u_{\star,m}-\Delta u_{\star,m}$ and for all complex $\tau$ satisfying $\text{Re}\,\tau>1$,  we have
$$\hat{y}(x,\tau)=\tau^\alpha \widehat{u_{\star,m}}(x,\tau)-\widehat{\Delta u_{\star,m}}(x,\tau)=0,\quad x\in\Omega.$$
Combining this with the uniqueness of the Laplace transform in time, we deduce that $\partial_t^\alpha u_{\star,m}-\Delta u_{\star,m}=0$ in $\Omega\times\Bbb R_+$. Finally, (1.5) implies that  for all complex $\tau$ satisfying $\text{Re}\,\tau>1$, we have 
$$\widehat{u_{\star,m}}(x,\tau)=\widehat{g_{\star,m}}(x,\tau),\quad x\in\partial\Omega$$
and applying again the uniqueness of Laplace transform in time we deduce that $u_{\star,m}=g_{\star,m}$ on $\partial\Omega\times \Bbb R_+$. Thus, $u_{\star,m}$ solves (1.1)  for $g=g_{\star,m}$. The uniqueness of  solutions $u$ of (1.1), with $g=g_{\star,m}$, satisfying $e^{- t}u_{\star,m}\in W^{1,\infty}(\Bbb R_+;H^2(\Omega))$ is the consequence of the uniqueness of the Laplace transform and the unique solvability of problem (1.5) for any complex $\tau$ satisfying $\text{Re}\,\tau>1$. Indeed, fix $u$ a solution of (1.1) with $g=0$ and satisfying $e^{- t}u\in W^{1,\infty}(\Bbb R_+;H^2(\Omega))$. Then, applying the Laplace transform in time to (1.1), we deduce that for any complex $\tau$ satisfying $\text{Re}\,\tau>1$,  $\hat{u}(\cdot,\tau)$ is well defined and it solves (1.5) with $g_{\star,m}=0$. Then the uniqueness of solutions of (1.5) implies that for all complex $\tau$ satisfying $\text{Re}\,\tau>1$,  $\hat{u}(\cdot,\tau)=0$ and the uniqueness of the Laplace transform in time implies that $u=0$.

Finally, since $e^{-t}u_{\star,m}\in L^{\infty}(\Bbb R_+;H^2(\Omega))$, for all $\tau>1$, we have
$$\begin{array}{ll}
\displaystyle
\int_0^{\infty}e^{-\tau t}\Vert\partial_{\nu}\,u_{\star,m}\Vert_{H^{\frac{1}{2}}(\partial\Omega)}\,dt
&
\displaystyle
\le
C\int_0^{\infty}e^{-(\tau-1)t}e^{-t}\Vert u_{\star,m}\Vert_{H^2(\Omega)}\,dt
\\
\\
\displaystyle
&
\displaystyle
\le C\Vert e^{-t}u_{\star,m}\Vert_{L^{\infty}(\Bbb R_+;H^2(\Omega))}<\infty.
\end{array}
$$
Therefore, for a.e. $x\in\partial\Omega$ and all $\tau>1$, we have
$$\begin{array}{ll}
\displaystyle
\widehat{\partial_{\nu}u_{\star,m}}\,(x,\tau)
&
\displaystyle
=\int_0^{\infty}e^{-\tau t}\,\partial_{\nu}\,u_{\star,m}(x,t)\,dt
\\
\\
\displaystyle
&
\displaystyle
=\partial_{\nu}
\left(\int_0^{\infty}e^{-\tau t}\,u_{\star,m}(x,t)\,dt\,\right)
\\
\\
\displaystyle
&
\displaystyle
=\partial_{\nu}\,\widehat{u_{\star,m}}\,(x,\tau).
\end{array}
$$
From this identity, we deduce (1.7).

\section{Proof of Theorem 1.1}

First we describe a basic system of inequalities which is a consequence of the expression (1.10) and the governing equations
of (1.3) and (1.9) for $w^0_{\star,m}$ and $w_{\star,m}$, respectively.

\proclaim{\noindent Lemma 3.1.}  
We have, for all $\tau> 1$
$$\begin{array}{c}
\displaystyle
I_{\star,m}(\tau)
\ge
\int_{\Omega}\frac{\tau^{\alpha_0}}{\tau^{\alpha(x)}}(\tau^{\alpha(x)}-\tau^{\alpha_0})(w_{\star,m}^0)^2\,dx
\end{array}
\tag {3.1}
$$
and
$$\displaystyle
I_{\star,m}(\tau)
\le\int_{\Omega}(\tau^{\alpha(x)}-\tau^{\alpha_0})(w_{\star,m}^0)^2\,dx.
\tag {3.2}
$$

\endproclaim

We omit to describe the proof since the idea of the derivation is well known in the framework of the enclosure method.
See \cite{Isize} and Proposition 4.1 in \cite{Iwall}.

It follows from (3.1) and (3.2) that

(0)  We have, for all $\tau> 1$
$$\displaystyle
\vert I_{\star,m}(\tau)\vert\le\tau^{\alpha_0}(\tau^{\Vert h\Vert_{L^{\infty}(D)}}+1)\,\int_D\,(w_{\star,m}^0(x))^2\,dx.
\tag {3.3}
$$
Note that the precise values of the power of $\tau$ is not important.

(i) if $\alpha>>\alpha_0$, then
$$\displaystyle
I_{\star,m}(\tau)\ge\tau^{\alpha_0}\int_{D}(\tau^{C\,\text{dist}\,(x,\partial D)^{\gamma}}-1)(w_{\star,m}^0(x))^2\,dx.
\tag {3.4}
$$

(ii) if $\alpha<<\alpha_0$, then
$$\displaystyle
I_{\star,m}(\tau)\le-\frac{\tau^{\alpha_0}}{\tau^{C\,\sup_{x\in D}\text{dist}(x,\partial D)^{\gamma}}}\int_{D}(\tau^{C\,\text{dist}\,(x,\partial D)^{\gamma}}-1)(w_{\star,m}^0(x))^2\,dx.
\tag {3.5}
$$

It is clear that $w^0_{\star,m}$ has the expression
$$\displaystyle
w^0_{\star,m}(x,\tau)
=\tau^{\alpha_0-1}v_{\star,m}(x;\beta)\vert_{\beta=\alpha_0},
\tag {3.6}
$$
where the function $v_{\star,m}(x;\beta)$ of $x\in\Bbb R^3$ with $\beta>0$ takes the form
$$
\displaystyle
v_{\star,m}(x;\beta)=
\left\{
\begin{array}{ll}
\displaystyle
\frac{1}{4\pi}\int_{B_{\eta}}(\eta^2-\vert y-p\vert^2)^m
\frac{e^{-\tau^{\frac{\beta}{2}}\vert x-y\vert}}{\vert x-y\vert}\,dy, & \text{if $\star=\text{ext}$,}
\\
\\
\displaystyle
\frac{1}{4\pi}\int_{B_{R_2}\setminus B_{R_1}}(R_2^2-\vert y-p\vert^2)^m(\vert y-p\vert^2-R_1^2)^m
\frac{e^{-\tau^{\frac{\beta}{2}}\vert x-y\vert}}{\vert x-y\vert}\,dy, & \text{if $\star=\text{int}$}
\end{array}
\right.
$$
and the point $p$ denotes the center of $B_{\eta}$ when $\star=\text{ext}$; the common center of $B_{R_1}$ and $B_{R_2}$
when $\star =\text{int}$.

\proclaim{\noindent Lemma 3.2.}  Let $\beta>0$, $\gamma\ge 0$, $C>0$ and $m\ge 0$ be an integer.
Then, there exist positive numbers $\tau_0>1$, $C_1$, $C_2$ and $\lambda\in\Bbb R$ such that, for all $\tau\ge\tau_0$
we have
$$\displaystyle
\tau^{\lambda}e^{2\tau^{\frac{\beta}{2}}\text{dist}\,(K_{\star},D)}
\int_{D}\,(\tau^{C\text{dist}\,(x,\partial D)^{\gamma}}-1)\,v_{\star,m}(x;\beta)^2\,dx\ge C_1
\tag {3.7}
$$
and
$$\displaystyle
\int_D\,v_{\star,m}(x;\beta)^2\,dx
\le C_2\tau^{-\beta(m+2)}e^{-2\tau^{\frac{\beta}{2}}\text{dist}\,(K_{\star},D)}.
\tag {3.8}
$$

\endproclaim

We postpone the proof of Lemma 3.2 to the next section.
We continue to prove Theorem 1.1.

First,  following (3.3), (3.6) and (3.8), we find
$$\displaystyle
\vert I_{\star,m}(\tau)\vert
\le\tau^{\lambda_1}\,e^{-2\tau^{\frac{\alpha_0}{2}}\text{dist}\,(K_{\star},D)},
\tag {3.9}
$$
where 
$$\displaystyle
\lambda_1=\alpha_0+\Vert h\Vert_{L^{\infty}(D)}+2(\alpha_0-1)-\beta(m+2).
$$
From this we obtain $\lim_{\tau\rightarrow\infty}e^{\frac{\alpha_0}{2}T}I_{\star,m}(\tau)=0$ for $T<2\text{dist}\,(K_{\star},D)$.

Next consider the case when $\alpha>>\alpha_0$.  From (3.4), (3.6) and (3.7)
one has
$$\displaystyle
I_{\star,m}(\tau)
\ge
C\tau^{\lambda_2}\,e^{-2\tau^{\frac{\alpha_0}{2}}\text{dist}\,(K_{\star},D)},
\tag {3.10}
$$
where $C$ is a positive constant independent of $\tau$ and
$$\displaystyle
\lambda_2=\alpha_0+2(\alpha_0-1)-\lambda.
\tag {3.11}
$$
A combination (3.9) and (3.10) yields (1.12) provided $\alpha>>\alpha_0$.
Besides, one has 
$$\displaystyle
e^{\tau^{\frac{\alpha_0}{2}}T}I_{\star,m}(\tau)\ge C_1\tau^{\lambda_2}e^{\tau^{\frac{\alpha_0}{2}}(T-2\,\text{dist}\,(K_{\star},D))}.
\tag {3.12}
$$
This yields $\lim_{\tau\rightarrow\infty}e^{\frac{\alpha_0}{2}T}I_{\star,m}(\tau)=\infty$ for $T>2\text{dist}\,(K_{\star},D)$
provided $\alpha>>\alpha_0$.

The proof in the case $\alpha<<\alpha_0$ can be done similarly as follows.
From (3.5), (3.6) and (3.7) we obtain
$$\displaystyle
I_{\star,m}(\tau)
\le
-C'\tau^{\lambda_3}\,e^{-2\tau^{\frac{\alpha_0}{2}}\text{dist}\,(K_{\star},D)},
\tag {3.13}
$$
where $C'$ is a positive constant independent of $\tau$ and 
$$\displaystyle
\lambda_3=\alpha_0-C\,\sup_{x\in D}\,\text{dist}(x,D)^{\gamma}+2(\alpha_0-1)-\lambda.
\tag {3.14}
$$
Note that $C$ above is the same one in the condition (A.II).
A combination (3.9) and (3.13) yields (1.12) provided $\alpha<<\alpha_0$.
Besides, one has 
$$\displaystyle
e^{\tau^{\frac{\alpha_0}{2}}T}I_{\star,m}(\tau)\le -C'\tau^{\lambda_3}e^{\tau^{\frac{\alpha_0}{2}}(T-2\,\text{dist}\,(K_{\star},D))}.
\tag {3.15}
$$
This yields $\lim_{\tau\rightarrow\infty}e^{\frac{\alpha_0}{2}T}I_{\star,m}(\tau)=-\infty$ for $T>2\text{dist}\,(K_{\star},D)$
provided $\alpha<<\alpha_0$.

This completes the proof of Theorem 1.1.

\section{Proof of Lemma 3.2}

Set $\tilde{\tau}=\tau^{\frac{\beta}{2}}$ and $v_{\star,m}(x;\beta)=v_{\star,m}(x)$.
By Proposition 3.1 in \cite{ISt}\footnote{Therein the case $m=0$ is excluded.  However, the proof still works also for the case.} 
we obtain the following:

\noindent
(i)  For $\vert x-p\vert>\eta$,  we have
$$\begin{array}{ll}
\displaystyle
v_{ext,m}(x)
&
\displaystyle
=\eta^{2(1+m)}
\cdot\frac{e^{-\tilde{\tau}\,\vert x-p\vert}}{\,\vert x-p\vert}\,
a_m(\tilde{\tau}),
\end{array}
\tag {4.1}
$$
where
$$\displaystyle
a_m(\tilde{\tau})
=\frac{1}{\tilde{\tau}}\int_0^1s(1-s^2)^m\sinh(\eta\tilde{\tau}\,s)\,ds.
$$

\noindent
(ii)  For $\vert x-p\vert<R_1$,  we get
$$\displaystyle
v_{int,m}(x)
=2
\cdot\frac{\sinh\,\tilde{\tau}\,\vert x-p\vert}{\vert x-p\vert}
\,b_m(\tilde{\tau}),
\tag {4.2}
$$
where
$$\displaystyle
b_m(\tilde{\tau})
=\frac{1}{\tilde{\tau}}\,\int_{R_1}^{R_2}\,s(R_2^2-s^2)^m(s^2-R_1^2)^m\,e^{-s\tilde{\tau}}\,ds.
$$
Besides, by Theorem 7.1 on p.81 in \cite{Ol}, we have,
as $\tilde{\tau}\rightarrow\infty$
$$
\displaystyle
a_m(\tilde{\tau})
\sim
\eta\,2^{m-1}m!\frac{e^{\tilde{\tau}\,\eta}}{(\tilde{\tau}\eta)^{m+2}}
\tag {4.3}
$$
and
$$
\\
\displaystyle
b_m(\tilde{\tau})
\sim
2^m\,m!R_1^{m+1}\,(R_2^2-R_1^2)^m
\frac{e^{-R_1\tilde{\tau}}}{\tilde{\tau}^{m+2}}.
\tag{4.4}
$$

\subsection{The case when $\star=ext$}
A combination of (4.1) and (4.3) gives, for all $x\in\Bbb R^3\setminus\overline B_{\eta}$
$$\displaystyle
v_{\text{ext},m}(x)^2\ge C_2^2\,\tau^{-\beta(m+2)}\frac{e^{-2\tau^{\frac{\beta}{2}}(\vert x-p\vert-\eta)}}{\vert x-p\vert^2}
\tag {4.5}
$$
and
$$\displaystyle
v_{\text{ext},m}(x)^2\le C_3^2\,\tau^{-\beta(m+2)}\frac{e^{-2\tau^{\frac{\beta}{2}}(\vert x-p\vert-\eta)}}{\vert x-p\vert^2},
\tag {4.6}
$$
where $C_2$ and $C_3$ are positive constants independent of $\tau$.

From (4.6) and the fact that $\inf_{x\in D}\vert x-p\vert-\eta=\text{dist}\,(\overline{B_{\eta}},D)$, we obtain (3.8) for $\star=\text{ext}$.

From (4.5) one gets
$$\begin{array}{l}
\,\,\,\,\,\,
\displaystyle
\tau^{\beta\,(m+2)}\,\int_{D}(\tau^{C\,\text{dist}\,(x,\partial D)^{\gamma}}-1)v_{\text{ext},m}(x)^2\,dx\\
\\
\displaystyle
\ge C_2^2\int_D\,(\tau^{C\,\text{dist}\,(x,\partial D)^{\gamma}}-1)
\frac{e^{-2\tau^{\frac{\beta}{2}}(\vert x-p\vert-\eta)}}{\vert x-p\vert^2}dx\\
\\
\displaystyle
\ge
\left(\frac{C_2}{\sup_{x\in D}\,\vert x-p\vert}\right)^2
\int_D\,(\tau^{C\,\text{dist}\,(x,\partial D)^{\gamma}}-1)\,e^{-2\tau^{\frac{\beta}{2}}(\vert x-p\vert-\eta)}dx.
\end{array}
\tag {4.7}
$$

First consider the case when $\gamma=0$.
By  Lemma A.2 in \cite{IThermo} under the assumption that $\partial D$ is $C^2$ we have
$$\displaystyle
(\tau^{\frac{\beta}{2}})^2e^{2\tau^{\frac{\beta}{2}}\,\text{dist}\,(D,\overline{B_{\eta}})}
\int_De^{-2\tau^{\frac{\beta}{2}}(\vert x-p\vert-\eta)}dx\ge C'.
$$
Thus, from (4.7) one gets
$$\displaystyle
\tau^{\beta\,(m+2)+\beta-C}e^{2\tau^{\frac{\beta}{2}}\,\text{dist}\,(D,\overline{B_{\eta}})}\int_{D}\,(\tau^{C}-1)\,v_{\text{ext},m}(x)^2\,dx
\ge \left(\frac{C_2}{\sup_{x\in D}\,\vert x-p\vert}\right)^2C'.
\tag {4.8}
$$

Next consider the case when $\gamma>0$.  We make a reduction to a simple geometry along the lines of the proof of Lemma A.1 in \cite{IThermo}.
Let $\tau\ge 1$.  Choose a point $q\in\partial D$ such that $\vert q-p\vert=d_{\partial D}(p)$.
Since $\partial D$ is $C^2$, one can find an open ball $B'$ with radius $\delta$ and centered at $q-\delta\nu_q$ such that
$B'\subset D$ and $\partial B'\cap\partial D=\{q\}$.  Then $\text{dist}\, (B',\overline{B_{\eta}})=\text{dist}\,(D,\overline{B_{\eta}})$
and $\text{dist}\,(x,\partial B')\le\text{dist}\,(x,\partial D)$ for all $x\in B'$.  
Thus, for all $x\in B'$ we have
$\tau^{C\text{dist}\,(x,\partial D)^{\gamma}}\ge\tau^{C\text{dist}\,(x,\partial B')^{\gamma}}\ge 1$.
Therefore, it suffices to prove (3.7) in the case when $D=B'$.

Set $d=d_{\partial D}(p)$.  Let $B''$ be the open ball with radius $d+\delta$ centered at $p$.
As described in the proof of Lemma A.1 in \cite{IThermo}, we have the global parametrization of
$B''\cap B'$:
$$\displaystyle
B''\cap B'
=\{\Upsilon(s,r,\theta)\,\vert\, 0<s<\delta, 0<r<(d+s)\sin\theta(s), \theta\in\,[0,2\pi[\},
$$
where\footnote{One can write
$$\displaystyle
\Upsilon(s,r,\theta)=p+(d+s)\omega,
$$
where
$$\displaystyle
\omega=\frac{1}{d+s}
\left(-\sqrt{(d+s)^2-r^2}\,\nu_q+r(\cos\theta \mbox{\boldmath $b$}+\sin\theta\mbox{\boldmath $c$}\right)\in S^2.
$$
}
$$\displaystyle
\Upsilon(s,r,\theta)
=p-\sqrt{(d+s)^2-r^2}\,\nu_q+r(\cos\theta \mbox{\boldmath $b$}+\sin\theta\mbox{\boldmath $c$});
$$
$\mbox{\boldmath $b$}$ and $\mbox{\boldmath $c$}$ are unit vectors chosen in such a way that
$\mbox{\boldmath $b$}\cdot\mbox{\boldmath $c$}=0$ and $\mbox{\boldmath $b$}\times\mbox{\boldmath $c$}=-\nu_q$;
$\theta(s)\in ]0,\,\frac{\pi}{2}[$ the unique solution of 
$$\displaystyle
\cos\theta
=\frac{(d+\delta)^2+(d+s)^2-\delta^2}{2(d+\delta)(d+s)}.
$$
We have
$$\displaystyle
\text{det}\,\Upsilon'(s,r,\theta)
=\frac{r(d+s)}
{\sqrt{(d+s)^2-r^2}}
$$
and
$$\begin{array}{ll}
\displaystyle
\text{dist}\,(x,\partial B')
&
\displaystyle
=\delta-\vert x-(p-(d+\delta)\nu_{q})\vert\\
\\
\displaystyle
&
\displaystyle
=\delta-\sqrt{\left(d+\delta-\sqrt{(d+s)^2-r^2}\,\right)^2+r^2},
\end{array}
$$
where $x=\Upsilon(s,r,\theta)$.  
The change of variables $x=\Upsilon(s,r,\theta)$ yields, for all $\tau\ge 1$
$$\begin{array}{l}
\,\,\,\,\,\,
\displaystyle
\int_{B'}\,(\tau^{C\,\text{dist}\,(x,\partial B')^{\gamma}}-1)\,e^{-2\tau^{\frac{\beta}{2}}(\vert x-p\vert-\eta)}dx
\\
\\
\displaystyle
\ge
\int_{B''\cap B'}\,(\tau^{C\,(\delta-\vert x-(p-(d+\delta)\nu_q)\vert)^{\gamma}}-1)\,e^{-2\tau^{\frac{\beta}{2}}(\vert x-p\vert-\eta)}dx
\\
\\
\displaystyle
=\int_0^{\delta}ds\int_0^{(d+s)\sin\theta(s)}dr
\int_0^{2\pi}d\theta
\frac{r(d+s)}
{\sqrt{(d+s)^2-r^2}}e^{-2\tau^{\frac{\beta}{2}}(d+s-\eta)}(\tau^{Cf(s,r)}-1)\\
\\
\displaystyle
=2\pi e^{-2\tau^{\frac{\beta}{2}}(d-\eta)}
\int_0^{\delta}ds\int_0^{(d+s)\sin\theta(s)}dr
\frac{r(d+s)}
{\sqrt{(d+s)^2-r^2}}e^{-2\tau^{\frac{\beta}{2}}s}(\tau^{Cf(s,r)}-1),
\end{array}
$$
where
$$\displaystyle
f(s,r)=\left\{\delta-\sqrt{\left(d+\delta-\sqrt{(d+s)^2-r^2}\right)^2+r^2}\,\right\}^{\gamma}.
$$
Making the change of a variable given by $r=(d+s)\sin\xi$, $\xi\in[0,\,\theta(s)]$, we have
$$\begin{array}{ll}
\displaystyle
I
&
\displaystyle
\equiv\int_0^{\delta}ds\int_0^{(d+s)\sin\theta(s)}dr
\frac{r(d+s)}
{\sqrt{(d+s)^2-r^2}}e^{-2\tau^{\frac{\beta}{2}}s}(\tau^{Cf(s,r)}-1)
\\
\\
\displaystyle
&
\displaystyle
=\int_0^{\delta}(d+s)^2e^{-2\tau^{\frac{\beta}{2}}s}ds
\int_0^{\theta(s)}
(\tau^{Cf(s,(d+s)\sin\xi)}-1)\sin\xi\,d\xi.
\end{array}
$$
Using 
$$\displaystyle
(d+\delta)^2+(d+s)^2=\delta^2+2(d+s)(d+\delta)\cos\theta(s),
$$
we have
$$\begin{array}{ll}
\displaystyle
f(s,(d+s)\sin\xi)
&
\displaystyle
=\left\{\delta-\sqrt{\left((d+\delta)-(d+s)\cos\xi\right)^2+(d+s)^2\sin^2\xi}\,\right\}^{\gamma}\\
\\
\displaystyle
&
\displaystyle
=\left\{\delta-\sqrt{(d+\delta)^2+(d+s)^2-2(d+\delta)(d+s)\cos\xi}\right\}^{\gamma}\\
\\
\displaystyle
&
\displaystyle
=\left\{\delta-\sqrt{\delta^2+2(d+\delta)(d+s)\cos\theta(s)-2(d+\delta)(d+s)\cos\xi}\right\}^{\gamma}\\
\\
\displaystyle
&
\displaystyle
=\left\{\delta-\sqrt{\delta^2-2(d+\delta)(d+s)(\cos\xi-\cos\theta(s))}\right\}^{\gamma}
\\
\\
\displaystyle
&
\displaystyle
=\left\{\frac{2(d+\delta)(d+s)(\cos\xi-\cos\theta(s))}
{\delta+\sqrt{\delta^2-2(d+\delta)(d+s)(\cos\xi-\cos\theta(s))}}
\right\}^{\gamma}
\\
\\
\displaystyle
&
\displaystyle
\ge
\left\{\frac{d+\delta}{\delta}(d+s)(\cos\xi-\cos\theta(s))
\right\}^{\gamma}
\\
\\
\displaystyle
&
\displaystyle
\ge
\left\{\frac{d(d+\delta)}{\delta}\right\}^{\gamma}
(\cos\xi-\cos\theta(s))^{\gamma}.
\end{array}
$$
This yields
$$\begin{array}{ll}
\displaystyle
\int_0^{\theta(s)}
(\tau^{Cf(s,(d+s)\sin\xi)}-1)\sin\xi\,d\xi
&
\displaystyle
\ge
\int_0^{\theta(s)}
(\tau^{C'(\cos\xi-\cos\theta(s))^{\gamma}}-1)\sin\xi\,d\xi\\
\\
\displaystyle
&
\displaystyle
=\int_0^{1-\cos\theta(s)}(\tau^{C'\sigma^{\gamma}}-1)d\sigma,
\end{array}
$$
where
$$\displaystyle
C'=C\left\{\frac{d(d+\delta)}{\delta}\right\}^{\gamma}.
$$
Here we have
$$\begin{array}{ll}
\displaystyle
1-\cos\theta(s)
&
\displaystyle
=\frac{s(2\delta-s)}{2(d+s)(d+\delta)}\\
\\
\displaystyle
&
\displaystyle
\ge
\frac{\delta s}{2(d+\delta)^2}.
\end{array}
$$
Thus one gets
$$
\displaystyle
\int_0^{\theta(s)}
(\tau^{Cf(s,(d+s)\sin\xi)}-1)\sin\xi\,d\xi
\ge
\int_0^{C''s}(\tau^{C'\sigma^{\gamma}}-1)d\sigma,
$$
where
$$\displaystyle
C''=\frac{\delta}{2(d+\delta)^2}.
$$
Therefore, we obtain
$$\displaystyle
I\ge d^2
\int_0^{\delta}e^{-2\tau^{\frac{\beta}{2}}s}ds
\int_0^{C''s}(\tau^{C'\sigma^{\gamma}}-1)d\sigma.
$$
Integrating by parts, we obtain
$$\begin{array}{l}
\,\,\,\,\,\,
\displaystyle
\int_0^{\delta}e^{-2\tau^{\frac{\beta}{2}}s}ds
\int_0^{C''s}(\tau^{C'\sigma^{\gamma}}-1)d\sigma\\
\\
\displaystyle
=-\frac{1}{2\tau^{\frac{\beta}{2}}}
\left(
[e^{-2\tau^{\frac{\beta}{2}}s}\int_0^{C''s}(\tau^{C'\sigma^{\gamma}}-1)d\sigma]_{s=0}^{s=\delta}
-C''\int_0^{\delta}
e^{-2\tau^{\frac{\beta}{2}}s}(\tau^{C_3s^{\gamma}}-1)ds\right)\\
\\
\displaystyle
=\frac{C''}{2\tau^{\frac{\beta}{2}}}
\int_0^{\delta}
e^{-2\tau^{\frac{\beta}{2}}s}(\tau^{C_3s^{\gamma}}-1)ds
-\frac{e^{-2\tau^{\frac{\beta}{2}}\delta}}{2\tau^{\frac{\beta}{2}}}
\int_0^{C''\delta}(\tau^{C'\sigma^{\gamma}}-1)d\sigma.
\end{array}
$$
Set
$$\displaystyle
C_3=C'(C'')^{\gamma}.
$$
Since we have
$$\displaystyle
\int_0^{C''\delta}(\tau^{C'\sigma^{\gamma}}-1)d\sigma
=O(\tau^{C_3\delta^{\gamma}}),
$$
one gets
$$\displaystyle
2\tau^{\frac{\beta}{2}}I
\ge
C''
\int_0^{\delta}
e^{-2\tau^{\frac{\beta}{2}}s}(\tau^{C_3s^{\gamma}}-1)ds
+O(\tau^{-\infty}).
$$
For simplicity, we write $C_3=C$.
We have
$$\begin{array}{ll}
\displaystyle
\,\,\,\,\,\,
\int_0^{\delta}
e^{-2\tau^{\frac{\beta}{2}}s}(\tau^{C_3s^{\gamma}}-1)ds & \displaystyle=-\frac{1}{2\tau^{\frac{\beta}{2}}}\int_0^{\delta} (e^{-2\tau^{\frac{\beta}{2}}s})'(\tau^{Cs^{\gamma}}-1)ds\\
\\
\displaystyle
&
\displaystyle
=\frac{C\gamma}{2}\tau^{-\frac{\beta}{2}}\log\tau\int_0^{\delta} e^{-2\tau^{\frac{\beta}{2}}s}\tau^{Cs^{\gamma}}s^{\gamma-1}ds
+O(\tau^{-\frac{\beta}{2}}e^{-2\tau^{\frac{\beta}{2}}\delta}(\tau^{C\delta^{\gamma}}-1))\\
\\
\displaystyle
&
\displaystyle
=\frac{C\gamma}{2}\tau^{-\frac{\beta}{2}}
\left(\log\tau\int_0^{\delta} e^{-2\tau^{\frac{\beta}{2}}s}\tau^{Cs^{\gamma}}s^{\gamma-1}ds
+O(\tau^{-\infty})\right)\\
\\
\displaystyle
&
\displaystyle
=\frac{C\gamma}{2}\tau^{-\frac{\beta}{2}}\log\tau
\left(K(\tau)
+O(\tau^{-\infty})\right),
\end{array}
$$
where
$$\displaystyle
K(\tau)\equiv\int_0^{\delta} e^{-2\tau^{\frac{\beta}{2}}s}\tau^{Cs^{\gamma}}s^{\gamma-1}ds.
$$
Let $\tau\ge 1$.  Since $\tau^{Cs^{\gamma}}\ge 1$, we have
$$
\displaystyle
K_p(\tau)
\ge
2^{-\gamma}\tau^{-\frac{\beta}{2}\gamma}\int_0^{2\tau^{\frac{\beta}{2}}\delta}
e^{-t}t^{\gamma-1}dt.
$$
This yields
$$\displaystyle
\liminf_{\tau\longrightarrow\infty}
\tau^{\frac{\beta}{2}\gamma}K(\tau)
\ge
2^{-\gamma}\int_0^{\infty}
e^{-t}t^{\gamma-1}dt
=\frac{\Gamma(\gamma)}{2^{\gamma}}.
$$
Therefore we obtain
$$\displaystyle
\liminf_{\tau\longrightarrow\infty}\tau^{\frac{\beta}{2}(\gamma+1)}(\log\tau)^{-1}\int_0^{\delta}
e^{-2\tau^{\frac{\beta}{2}}}(\tau^{C_3s^{\gamma}}-1)ds
\ge
\frac{C_3\Gamma(\gamma+1)}{2^{\gamma+1}},
$$
Thus, we obtain
$$\displaystyle
\liminf_{\tau\longrightarrow\infty}2\tau^{\frac{\beta}{2}(\gamma+2)}(\log\tau)^{-1}I
\ge \frac{C''C_3\Gamma(\gamma+1)}{2^{\gamma+1}}.
$$
Summing up, from this together with (4.8) we see that the $\lambda$ on (3.7) should be 
$$\displaystyle
\lambda=
\left\{
\begin{array}{ll}
\beta(m+3)-C & \text{if $\gamma=0$},\\
\\
\displaystyle
\beta(m+2)+\frac{\beta}{2}(\gamma+2)+\epsilon & \text{if $0<\gamma$,}
\end{array}
\right.
$$
where $\epsilon$ is an arbitrary positive number.

\subsection{The case when $\star=int$}

A combination of (4.2) and (4.4) gives, for all $x\in B_{R_1}$
$$\displaystyle
v_{\text{int},m}(x)^2\ge C_3^2\,\tilde{\tau}^{-2(m+2)}\left(\frac{\sinh\,\tilde{\tau}\,\vert x-p\vert}{\vert x-p\vert}\right)^2
e^{-2R_1\tilde{\tau}}
\tag {4.9}
$$
and
$$\displaystyle
v_{\text{int},m}(x)^2\le C_4^2\,\tilde{\tau}^{-2(m+2)}\left(\frac{\sinh\,\tilde{\tau}\,\vert x-p\vert}{\vert x-p\vert}\right)^2
e^{-2R_1\tilde{\tau}},
\tag {4.10}
$$
where $C_3$ and $C_4$ are positive constants independent of $\tau$.

From (4.10) and $R_1-R_D(p)=\text{dist}\,(\overline{B_{R_2}}\setminus B_{R_1}, D)$ we obtain (2.16) for $\star=\text{int}$.

We make a reduction to a simple geometry along the lines of the proof of Lemma 4.3 in \cite{IE06}.
Let $\tau\ge 1$.  The $D$ is contained in the open ball $B_{R_D(p)}$ centered at $p$ with radius $R_D(p)$ and $B_{R_D(p)}
\subset B_{R_1}$.
Choose a point $q\in\partial D$ such that $\vert q-p\vert=R_{D}(p)$.
Since $\partial D$ is $C^2$, one can find an open ball $B'$ with radius $\delta<\frac{\rho}{2}$ and centered at $q-\delta\nu_q$ such that
$B'\subset D$ and $\partial B'\cap\partial D=\{q\}$.  
Then $\text{dist}\,(\overline{B_{R_2}}\setminus B_{R_1},B')=R_1-R_D(p)=\text{dist}\,(\overline{B_{R_2}}\setminus B_{R_1},D)$
and $\text{dist}\,(x,\partial B')\le\text{dist}\,(x,\partial D)$ for all $x\in B'$.  
Thus, for all $x\in B'$ we have
$\tau^{C\text{dist}\,(x,\partial D)^{\gamma}}\ge\tau^{C\text{dist}\,(x,\partial B')^{\gamma}}\ge 1$.
Therefore it suffices to prove (3.7) in the case when $D=B'$.  Up to this point, it is the same as above.

Set $\rho=R_D(p)$ and let $0<\delta'<\delta$ and $B''$ be the open ball with radius $\rho-\delta'$ centered at $p$.
We make use of the parametrization of the set $B'\setminus\overline B''$ which is essentially same as that used in \cite{IE06}:
$$\displaystyle
B''\setminus\overline B'
=\cup_{0<s<\delta'}\left\{p+(\rho-s)\omega\,\vert\,\omega\in S(s)\right\},
$$
where
$$\displaystyle
S(s)=\left\{\omega\in S^2\,\vert\,\omega\cdot\nu_q>\cos\theta(s)
\right\}
$$
and the $\theta(s)\in \,]0,\,\frac{\pi}{2}[$ is the unique solution of the equation
$$\displaystyle
\cos\theta=\frac{(\rho-\delta')^2+(\rho-s)^2-\delta'^2}
{2(\rho-\delta')(\rho-s)}.
$$

Choose two linearly independent unit vectors $\mbox{\boldmath $b$}$ and $\mbox{\boldmath $c$}$
in such a way that $\mbox{\boldmath $b$}\cdot\mbox{\boldmath $c$}=0$ and $\mbox{\boldmath $b$}\times\mbox{\boldmath $c$}
=\nu_q$.  Then we have the expression
$$\displaystyle
B'\setminus\overline{B''}
=
\left\{
\Upsilon(s,r,\theta)\,\vert\,
0<s<\delta', 
0\le r<(\rho-s)\sin\theta(s),
0\le\theta<2\pi
\right\},
$$
where
$$
\displaystyle
\Upsilon\,(s,r,\theta)=p+\sqrt{(\rho-s)^2-r^2}\,\nu_q+r(\cos\theta\mbox{\boldmath $b$}+\sin\theta\mbox{\boldmath $c$}).
$$
We have
$$\displaystyle
\text{det}\,\Upsilon'(s,r,\theta)
=-\frac{r(\rho-s)}
{\sqrt{(\rho-s)^2-r^2}}
$$
and
$$\begin{array}{ll}
\displaystyle
\text{dist}\,(x,\partial B')
&
\displaystyle
=\delta'-\vert x-(p+(\rho-\delta')\nu_{q})\vert\\
\\
\displaystyle
&
\displaystyle
=\delta'-\sqrt{\left(
\rho-\delta'-\sqrt{(\rho-s)^2-r^2}\,\right)^2+r^2},
\end{array}
$$
where $x=\Upsilon(s,r,\theta)$.

The change of variables $x=\Upsilon(s,r,\theta)$ yields, for all $\tau\ge 1$
$$\begin{array}{l}
\,\,\,\,\,\,
\displaystyle
\int_{B'}\,(\tau^{C\,\text{dist}\,(x,\partial B')^{\gamma}}-1)\,
\left(\frac{\sinh\,\tilde{\tau}\,\vert x-p\vert}{\vert x-p\vert}\right)^2
e^{-2R_1\tilde{\tau}}\,dx
\\
\\
\displaystyle
\ge
\int_{B'\setminus B''}\,(\tau^{C\,(\delta'-\vert x-(p+(\rho-\delta')\nu_q)\vert)^{\gamma}}-1)\,\left(\frac{\sinh\,\tilde{\tau}\,\vert x-p\vert}{\vert x-p\vert}\right)^2
e^{-2R_1\tilde{\tau}}\,dx
\\
\\
\displaystyle
=\int_0^{\delta'}ds\int_0^{(\rho-s)\sin\theta(s)}dr
\int_0^{2\pi}d\theta
\frac{r(\rho-s)}
{\sqrt{(\rho-s)^2-r^2}}
\left(\frac{\sinh\,\tilde{\tau}\,(\rho-s)}{\rho-s}\right)^2
e^{-2R_1\tilde{\tau}}
(\tau^{Cf(s,r)}-1)\\
\\
\displaystyle
=2\pi e^{-2(R_1-\rho)\tilde{\tau}}
\int_0^{\delta'}ds\int_0^{(\rho-s)\sin\theta(s)}dr
\frac{r}
{\sqrt{(\rho-s)^2-r^2}}
\frac{e^{-2\rho\tilde{\tau}}\sinh^2\,\tilde{\tau}\,(\rho-s)}{\rho-s}
(\tau^{Cf(s,r)}-1),
\end{array}
\tag {4.11}
$$
where
$$\displaystyle
f(s,r)=\left\{\delta'-\sqrt{\left(
\rho-\delta'-\sqrt{(\rho-s)^2-r^2}\,\right)^2+r^2}\,
\right\}^{\gamma}.
$$
Making the change of variable
$$\begin{array}{ll}
\displaystyle
r=(\rho-s)\,\sin\xi,
& 
0<\xi<\theta(s),
\end{array}
$$
we have
$$\begin{array}{ll}
\displaystyle
I'
&
\displaystyle
\equiv\int_0^{\delta'}ds\int_0^{(\rho-s)\sin\theta(s)}dr
\frac{r}
{\sqrt{(\rho-s)^2-r^2}}
\frac{e^{-2\rho\tilde{\tau}}\sinh^2\,\tilde{\tau}\,(\rho-s)}{\rho-s}
(\tau^{Cf(s,r)}-1)
\\
\\
\displaystyle
&
\displaystyle
=\int_0^{\delta'}ds
e^{-2\rho\tilde{\tau}}\sinh^2\,\tilde{\tau}\,(\rho-s)
\int_0^{\theta(s)}(\tau^{Cf(s,(\rho-s)\sin\xi)}-1)\sin\xi\,d\xi.
\end{array}
$$
Here using 
$$\displaystyle
(\rho-\delta')^2+(\rho-s)^2=\delta'^2+2(\rho-s)(\rho-\delta')\cos\theta(s),
$$
we have
$$\begin{array}{ll}
\displaystyle
f(s,(\rho-s)\sin\xi)
&
\displaystyle
=\left\{\delta'-\sqrt{\left((\rho-\delta')-(\rho-s)\cos\xi\right)^2+(\rho-s)^2\sin^2\xi}\,\right\}^{\gamma}\\
\\
\displaystyle
&
\displaystyle
=\left\{\delta'-\sqrt{(\rho-\delta')^2+(\rho-s)^2-2(\rho-\delta')(\rho-s)\cos\xi}\right\}^{\gamma}\\
\\
\displaystyle
&
\displaystyle
=\left\{\delta'-\sqrt{\delta'^2+2(\rho-\delta')(\rho-s)\cos\theta(s)-2(\rho-\delta')(\rho-s)\cos\xi}\right\}^{\gamma}\\
\\
\displaystyle
&
\displaystyle
=\left\{\delta'-\sqrt{\delta'^2-2(\rho-\delta')(\rho-s)(\cos\xi-\cos\theta(s))}\right\}^{\gamma}
\\
\\
\displaystyle
&
\displaystyle
=\left\{\frac{2(\rho-\delta')(\rho-s)(\cos\xi-\cos\theta(s))}
{\delta'+\sqrt{\delta'^2-2(\rho-\delta')(\rho-s)(\cos\xi-\cos\theta(s))}}
\right\}^{\gamma}
\\
\\
\displaystyle
&
\displaystyle
\ge
\left\{\frac{\rho-\delta'}{\delta'}(\rho-s)(\cos\xi-\cos\theta(s))
\right\}^{\gamma}
\\
\\
\displaystyle
&
\displaystyle
\ge
\left\{\frac{(\rho-\delta')^2}{\delta'}\right\}^{\gamma}
(\cos\xi-\cos\theta(s))^{\gamma}.
\end{array}
$$
This yields
$$\begin{array}{ll}
\displaystyle
\int_0^{\theta(s)}
(\tau^{Cf(s,(\rho-s)\sin\xi)}-1)\sin\xi\,d\xi
&
\displaystyle
\ge
\int_0^{\theta(s)}
(\tau^{C'(\cos\xi-\cos\theta(s))^{\gamma}}-1)\sin\xi\,d\xi\\
\\
\displaystyle
&
\displaystyle
=\int_0^{1-\cos\theta(s)}(\tau^{C''\sigma^{\gamma}}-1)d\sigma,
\end{array}
$$
where
$$\displaystyle
C''=C\left\{\frac{(\rho-\delta')^2}{\delta'}\right\}^{\gamma}.
$$
Here we have
$$\begin{array}{ll}
\displaystyle
1-\cos\theta(s)
&
\displaystyle
=\frac{s(2\delta'-s)}{2(\rho-s)(\rho-\delta')}\\
\\
\displaystyle
&
\displaystyle
\ge
\frac{\delta' s}{2\rho(\rho-\delta')}.
\end{array}
$$
Thus one gets
$$
\displaystyle
\int_0^{\theta(s)}
(\tau^{Cf(s,(\rho-s)\sin\xi)}-1)\sin\xi\,d\xi
\ge
\int_0^{C'''s}(\tau^{C''\sigma^{\gamma}}-1)d\sigma,
$$
where
$$\displaystyle
C'''=\frac{\delta}{2\rho(\rho-\delta')}.
$$
Therefore, we obtain
$$\displaystyle
I'\ge
e^{-2\rho\tilde{\tau}}
\int_0^{\delta'}ds
\sinh^2\,\tilde{\tau}\,(\rho-s)
\int_0^{C'''s}(\tau^{C''\sigma^{\gamma}}-1)d\sigma.
\tag {4.12}
$$
Since we have
$$\displaystyle
\int\,\sinh^2\,x\, dx
=\frac{e^{2x}-e^{-2x}}{8}-\frac{x}{2},
$$
integration by parts yields
$$
\begin{array}{l}
\displaystyle
\,\,\,\,\,\,
e^{-2\rho\tilde{\tau}}
\int_0^{\delta'}ds
\sinh^2\,\tilde{\tau}\,(\rho-s)
\int_0^{C'''s}(\tau^{C''\sigma^{\gamma}}-1)d\sigma
\\
\\
\displaystyle
=
-\tilde{\tau}^{-1}e^{-2\rho\tilde{\tau}}
\left\{
\left[
\left(\frac{e^{2\tilde{\tau}(\rho-s)}-e^{-2\tilde{\tau}(\rho-s)}}{8}-\frac{\tilde{\tau}(\rho-s)}{2}\right)
\int_0^{C'''s}(\tau^{C''\sigma^{\gamma}}-1)d\sigma\,\right]_{s=0}^{s=\delta'}\right.
\\
\\
\displaystyle
\,\,\,
\left.
-\int_0^{\delta'}\,\left(\frac{e^{2\tilde{\tau}(\rho-s)}-e^{-2\tilde{\tau}(\rho-s)}}{8}-\frac{\tilde{\tau}(\rho-s)}{2}\right)
C^{'''}(\tau^{C''s^{\gamma}}-1)\,ds
\right\}
\\
\\
\displaystyle
=
-\tilde{\tau}^{-1}e^{-2\rho\tilde{\tau}}
\left(\frac{e^{2\tilde{\tau}(\rho-\delta')}-e^{-2\tilde{\tau}(\rho-\delta')}}{8}-\frac{\tilde{\tau}(\rho-\delta')}{2}\right)
\int_0^{C'''\delta'}(\tau^{C''\sigma^{\gamma}}-1)d\sigma
\\
\\
\displaystyle
\,\,\,
+\tilde{\tau}^{-1}e^{-2\rho\tilde{\tau}}C^{'''}
\int_0^{\delta'}\,\left(\frac{e^{2\tilde{\tau}(\rho-s)}-e^{-2\tilde{\tau}(\rho-s)}}{8}-\frac{\tilde{\tau}(\rho-s)}{2}\right)
(\tau^{C''s^{\gamma}}-1)\,ds\\
\\
\displaystyle
=
O(\tilde{\tau}^{-1}e^{-2\tilde{\tau}\delta'}+\tilde{\tau}^{-1}e^{-2\tilde{\tau}(2\rho-\delta')}
+e^{-2\rho\tilde{\tau}})
\int_0^{C'''\delta'}(\tau^{C''\sigma^{\gamma}}-1)d\sigma
\\
\\
\displaystyle
\,\,\,
+\tilde{\tau}^{-1}\frac{C^{'''}}{8}
\int_0^{\delta'}\,e^{-2\tilde{\tau}s}(\tau^{C''s^{\gamma}}-1)\,ds
+O(\tilde{\tau}^{-1}e^{-2\tilde{\tau}(2\rho-\delta')}+e^{-2\tilde{\tau}\rho})\int_0^{\delta'}(\tau^{C''s^{\gamma}}-1)ds
\end{array}
\tag {4.13}
$$
Here we have
$$\displaystyle
\int_0^{C'''\delta'}(\tau^{C''\sigma^{\gamma}}-1)d\sigma
=O(\tau^{C_4})
$$
and
$$\displaystyle
\int_0^{\delta'}(\tau^{C''s^{\gamma}}-1)ds
=O(\tau^{C_5}),
$$
where
$$\begin{array}{ll}
\displaystyle
C_4=C^{''}(C^{'''}\delta')^{\gamma}, 
&
\displaystyle
C_5=C^{''}(\delta')^{\gamma}.
\end{array}
$$
Therefore from this together with (4.12) and (4.13) we obtain
$$\displaystyle
I'\ge
\tilde{\tau}^{-1}\frac{C^{'''}}{8}
\int_0^{\delta'}\,e^{-2\tilde{\tau}s}(\tau^{C''s^{\gamma}}-1)\,ds
+O(\tau^{-\infty}).
\tag {4.14}
$$
Consider the case when $\gamma=0$.  From (4.14) we have
$$\displaystyle
I'\ge
\tilde{\tau}^{-1}\frac{C^{'''}}{8}(\tau^{C''}-1)
\int_0^{\delta'}\,e^{-2\tilde{\tau}s}\,ds
+O(\tau^{-\infty})
$$
and thus
$$\displaystyle
\liminf_{\tau\rightarrow\infty}\tau^{\beta-C''}I'>0.
\tag {4.15}
$$
Next let $\gamma>0$.  It follows from the case when $\star=\text{ext}$ and $\gamma>0$ one has
$$\displaystyle
\liminf_{\tau\longrightarrow\infty}\tau^{\frac{\beta}{2}(\gamma+1)}(\log\tau)^{-1}\int_0^{\delta'}
e^{-2\tilde{\tau}s}(\tau^{C^{''}s^{\gamma}}-1)ds
\ge
\frac{C_3\Gamma(\gamma+1)}{2^{\gamma+1}},
$$
Thus from (4.14) one gets
$$\displaystyle
\liminf_{\tau\rightarrow\infty}8\tau^{\frac{\beta}{2}(\gamma+2)}(\log\tau)^{-1}I'\ge
\frac{C^{'''}C_3\Gamma(\gamma+1)}{2^{\gamma+1}}.
\tag {4.16}
$$
Noting $R_1-\rho=\text{dist}\,(K_{\text{int}},D)$, from (4.9), (4.11), (4.15) and (4.16) we see that the $\lambda$ on (3.7) should be 
$$\displaystyle
\lambda=
\left\{
\begin{array}{ll}
\beta(m+3)-C'' & \text{if $\gamma=0$},\\
\\
\displaystyle
\beta(m+2)+\frac{\beta}{2}(\gamma+2)+\epsilon & \text{if $0<\gamma$,}
\end{array}
\right.
$$
where $\epsilon$ is an arbitrary positive number.

\section{Some additional remark}

(1)  It seems that the growth order of the absolute value of the function $e^{\frac{\alpha_0}{2}T}\,I_{\star,m}(\tau)$ in Theorem 1.1
as $\tau\rightarrow\infty$ for $T>2\text{dist}\,(K_{\star},D)$ becomes worse if $\gamma$ or $m$ is large.
See (3.12) with (3.11), (3.15) with (3.14) and $\lambda$ with $\beta=\alpha_0$ in the end of Subsection 4.1 in the case $\star=\text{ext}$
and Subsection 4.2 in the case $\star=\text{int}$.

(2)  Replace $g_{\star,m}$ in Proposition 1.1 with $g_{\star,m_1,m_2}$ given by
$$\begin{array}{lll}
\displaystyle
g_{\star,m_1,m_2}(x,t)=\frac{e^t}{2\pi}\,
\int_{-\infty}^{\infty}
e^{its}(1+is)^{-5}w_{\star,m_1,m_2}^0(x,1+is)ds, & x\in\partial\Omega, & t\in[0,\,\infty[,
\end{array}
$$
where the function $w^0_{\star,m_1,m_2}\in H^2(\Bbb R^3)$ is the unique solution of the equation
$$\begin{array}{ll}
\displaystyle
(\Delta-\tau^{\alpha_0}\,)w_{\star,m_1,m_2}^0+\tau^{\alpha_0-1}\,\Phi_{\star}(x)=0, & \displaystyle x\in\Bbb R^3,
\end{array}
$$
and $\Phi_{\star}$ an arbitrary fixed measurable function of $x\in\Bbb R^3$ such that
$$\begin{array}{ll}
\displaystyle
C^{-1}\Psi_{\star,m_1}(x)\le\Phi_{\star}(x)\le C\,\Psi_{\star,m_2}(x) & \text{a.e. $x\in\Bbb R^3$}
\end{array}
\tag {5.1}
$$
with some $m_1,m_2=0,1,\cdots$ and $C>0$.

It is easy to see that the proof of Proposition 1.1 still works also for this case
since the concrete form of $\Phi_{\star}$ is not used therein.

Then, one can define a new indication function by trivial replacements of $w^0_{\star,m}$ and $u_{\star,m}$ on (1.8).
We see that Lemma 3.1 where the old indicator function replaced with this new one,
is valid by replacing $w^0_{\star,m}$ on (3.1) and (3.2) with $w^0_{\star,m_1,m_2}$.
Besides, by virtue of assumption (5.1), for $\tau>1$, we have
$$\begin{array}{ll}
\displaystyle
w^0_{\star,m_1}(x)\le w^0_{\star,m_1,m_2}(x)\le w^0_{\star, m_2}(x) & \text{a.e.$x\in\Bbb R^3$.}
\end{array}
$$ 
Thus, the estimates for the new indicator function corresponding to (3.3), (3.4) and (3.5)
are valid under the replacements: $m$ of $W^0_{\star,m}$ on (3.3) with $m_2$
and (3.4) and (3.5) with $m_1$.
Hereafter the proof of Theorem 1.1 works also for the new indicator function and one gets 
the corresponding result to Theorem 1.1.

$$\quad$$

\centerline{{\bf Acknowledgments}}

MI was partially supported by the Grant-in-Aid for
Scientific Research (C)(No. 17K05331) and (B)(No. 18H01126) of Japan  Society for the
Promotion of Science.  The work of YK was partially supported by  the French National Research Agency ANR (project MultiOnde) grant ANR-17-CE40-0029.

\end{document}